\documentclass[12pt, reqno]{amsart}
\usepackage{latexsym}
\usepackage{amsfonts}
\usepackage{amssymb}
\usepackage{graphicx}
\usepackage{graphicx}
\usepackage{color}

\newtheorem{theorem}{Theorem}[section]

\newtheorem{cor}[theorem]{Corollary}
\newtheorem{rem}[theorem]{Remark}
\numberwithin{equation}{section}

\newcommand{\bC}{{\mathbb{C}}}

\newcommand{\bM}{{\mathbb{M}}}

\newcommand{\bR}{{\mathbb{R}}}

\newcommand{\spn}{\operatorname{span}}

\newcommand{\tr}{\operatorname{tr}}
\newcommand{\supp}{\operatorname{supp}}

\newcommand{\CC}{\mathbb{C}}

\newcommand{\EE}{\mathcal{E}}

\begin{document}
\title[Functional Equations and Fourier Analysis]
{Functional Equations and Fourier Analysis}
\author[D. Yang]
{Dilian Yang}
\address{DILIAN YANG, Department of Mathematics $\&$ Statistics,
University of Windsor, Windsor, ON N9B 3P4,
CANADA}
\email{dyang@uwindsor.ca}

\begin{abstract}
By exploring the relations among functional equations, harmonic analysis and representation theory,
we give a unified and very accessible approach to solve three important functional equations --
the d'Alembert equation, the Wilson equation, and the d'Alembert long equation,
on compact groups.
\end{abstract}

\thanks{The research was partially supported by an NSERC Discovery grant.}

\date{}
\maketitle

Recently,  three important equations --
the d'Alembert equation, the Wilson equation, and the d'Alembert long equation
have been attracting a great deal of attention.
See \cite{AnYang1, AnYang2, ABE, Cho, Dav1, Dav2, Stet4, Stet5, Yang06, Yangthesis} and the references therein.
It turns out that their solutions have very nice structures.
In particular, their solutions on compact groups were obtained as consequences
of the main results of \cite{AnYang2}, where
a much more general class of functional equations was studied.

Because of the increasing importance of these three equations,
it is worthwhile to give a much more transparent approach. This is the first main purpose of this
short note:
We shall give a unified, transparent, and very accessible approach to solve these
equations on compact groups.

The second main purpose is, as in \cite{AnYang1, AnYang2, Yang06, Yangthesis},
to explore the relations among different areas whenever possible, such as the areas of
functional equations, harmonic analysis and representation theory considered here.

 \section{Preliminaries}\label{S:Pre}

In this section, we set up some notation and conventions,
briefly review some fundamental facts in Fourier analysis
which will be used later, and
introduce the functional equations we shall be concerned with.

\subsection{Fourier analysis}
\label{SS:fa}

Let $G$ be a compact group with the normalized Haar measure $dx$. Let $\hat{G}$
stand for the set of equivalence classes of irreducible unitary representations of $G$.
For $[\pi]\in\hat{G}$,
the notation $d_\pi$ denotes the dimension of the representation space of $\pi$.
For $f\in L^2(G)$, the \textit{Fourier transform} of $f$ is defined by
$$
\hat{f}(\pi)=d_\pi\int_Gf(x)\pi(x)^{-1}dx\in\bM_{d_\pi}(\CC)
\quad \text{for all} \quad [\pi]\in\hat{G},
$$
where
$\bM_n(\CC)$ is the space of all $n\times n$ complex matrices.

As usual, the left and right regular representations of
$G$ in $L^2(G)$ are defined by
$$(L_yf)(x)=f(y^{-1}x), \quad (R_yf)(x)=f(xy),$$ respectively, where
$f\in L^2(G)$ and $x,y\in G$.
A crucial property of the Fourier transform is that it converts the
regular representations of $G$ into matrix multiplications.

The following facts will be useful later.
\begin{itemize}
\item[{(i)}] The
\textit{Fourier inversion formula} is given by
$$f(x)=\sum_{[\pi]\in\hat{G}}\tr(\hat{f}(\pi)\pi(x))
\quad \text{for all} \quad x\in G.$$
\item[{(ii)}] The following identities hold:
\begin{align*}
(L_yf)\hat{}\,(\pi)=\hat{f}(\pi)\pi(y)^{-1}, \quad
(R_yf)\hat{}\,(\pi)=\pi(y)\hat{f}(\pi)
\end{align*}
for all $y\in G$, and $\pi\in\hat G$.
\end{itemize}

For more information about the topics of this subsection, refer to \cite[Chapter 5]{Forlland}.

\subsection{Functional equations}
Let $G$ be a compact group, and $f$, $g$ be complex functions on $G$.
In this note, we study the following functional equations
\begin{align}
\label{E:D'Alembert}
f(xy)+f(xy^{-1})&=2f(x)f(y), \\
\label{E:Wilson}
f(xy)+f(xy^{-1})&=2f(x)g(y), \\
\label{E:long}
f(xy)+f(yx)+f(xy^{-1})+f(y^{-1}x)&=4f(x)f(y),
\end{align}
which are
called the \textit{d'Alembert equation},
the \textit{Wilson equation}, and the \textit{d'Alembert long equation}, respectively.
It is known that those equations are closely related to each other. For example,
Eq.~\eqref{E:D'Alembert} is a special case of Eq.~\eqref{E:Wilson};
Eq.~\eqref{E:long} becomes Eq.~\eqref{E:D'Alembert} if $f$ in Eq.~\eqref{E:long} is central
(i.e., $f(xy)=f(yx)$ for all $x,y\in G$);
and if $(f,g)$ with $f\ne 0$ satisfies Eq.~\eqref{E:Wilson} then $g$ also satisfies Eq.~\eqref{E:long}.

A lot of effort has been put on solving those equations recently. It turns out that
their solutions have very nice and interesting structures. See, e.g., \cite{AnYang1, Dav1, Dav2, Stet5} for
more details.

\smallskip
Throughout the rest of this note, the group $G$ is always assumed to be compact.
By solutions (resp. representations), we always mean continuous solutions
(resp. continuous representations).

\section{Small dimension lemma}
\label{S:sdl}

We will prove a very useful lemma in this section, which may be of independent interest.
The lemma, roughly speaking, says that for an \textit{irreducible} representation
$\pi$ of a compact group $G$, if the operators $\pi(x)+\pi(x)^{-1}$ for all $x\in G$
have a common nonzero eigenvector, then the dimension $d_\pi$ of $\pi$ has to be rather small:
$d_\pi\le 2$.
For obvious reasons, we call this lemma \textit{Small Dimension Lemma}.
In the next section, we will apply it to give a unified approach to solve
Eqs.~\eqref{E:D'Alembert}, \eqref{E:Wilson}, \eqref{E:long}.

\medskip
\noindent
\textbf{Small Dimension Lemma.}
\textit{
Let $G$ be a compact group, and $\pi:G\rightarrow U(n)$ an
irreducible representation of $G$. Suppose that there exists a nonzero
vector $v\in\CC^n$ such that
\begin{align}
\label{lemcon}
(\pi(x)+\pi(x)^{-1})v\in\CC v \quad \mbox{for all}\quad  x\in G.
\end{align}
Then either $n=1$ and $\pi$ is a unitary group character,
or $n=2$ and $\pi(G)\subseteq SU(2)$.
}

\medskip

\begin{proof}
For $x\in G$, let $\pi_{ij}(x)=\langle\pi(x)e_j,e_i\rangle$
denote the $(i,j)$-th entry of
$\pi(x)$. Consider the subspaces $\EE_i$ of $L^2(G)$:
$$
\EE_i=\spn\{\pi_{ij}\mid j=1,\ldots,n\}, \quad i=1,\ldots,n.
$$
Let $R|_{\EE_i}$ be the restriction of the right regular
representation $R$ of $G$ to $\EE_i$. Then $R|_{\EE_i}$ is equivalent
to $\pi$. Indeed, let $U:\EE_i\to \CC^n$ be the unitary operator defined via $U(\pi_{ij})=e_j$;
then one can easily check that $U^*\pi(\cdot) U=R|_{\EE_i}(\cdot)$.
In particular, $R|_{\EE_i}$ is irreducible as $\pi$ is irreducible.
Thus
$$
\EE_i=\spn\{R_x\pi_{ij}\mid x\in G\}\quad \text{for any}\quad j=1,\ldots,n.
$$

After similarity, we may and do assume that the nonzero vector $v$
in \eqref{lemcon} is given by
$v=(1,0,\ldots,0)^t$, where $t$ denotes transpose.
From the condition \eqref{lemcon}, it follows that
\begin{eqnarray}
\label{c:1}
\pi_{i1}=-\overline{\pi_{1i}}\quad \text{for all}\quad i=2,\ldots,n.
\end{eqnarray}
Here $\overline{\pi_{1i}}$ denotes the complex conjugate of $\pi_{1i}$.
Hence we have
\begin{alignat*}{2}
\EE_i&=\spn\{R_x\pi_{i1}\mid x\in G\}
      =\spn\{R_x\overline{\pi_{1i}}\mid x\in G\}\\
     &=\spn\{\overline{R_x\pi_{1i}}\mid x\in G\}=\overline{\EE}_1
\end{alignat*}
for each  $i=2, \ldots, n$.
If $n>2$, then $\EE_2=\EE_3=\cdots=\EE_n$.  This is impossible by Schur's
orthogonality relations (cf. \cite{Forlland}). Thus $n\leq 2$.

In the case of $n=1$, $\pi$ is of course a unitary group character.

In the case of $n=2$, it follows from $\pi_{21}=-\overline{\pi_{12}}$ in \eqref{c:1} that
$\pi(x)=\left(\begin{matrix}a&b\\  -\bar{b}&c\end{matrix}\right)$.
As $\pi(x)$ is a unitary operator, some simple calculations give either
$\pi(x)=\left(\begin{matrix}a&b\\ -\overline{b}&\overline{a}\end{matrix}\right)\in SU(2)$
with $b\ne 0$,
or
$\pi(x)=\left(\begin{matrix}a&0\\ 0&c\end{matrix}\right)$.
But the product of any such two elements should also be one of the forms.
This forces that either $\pi(x)\in SU(2)$ for all $x\in G$,  or
$\pi(x)\in \left\{\left(\begin{matrix}a&0\\ 0&c\end{matrix}\right)\in U(2)\right\}$ for all
$x\in G$. Since $\pi$ is irreducible, we necessarily have $\pi(G)\subseteq SU(2)$.
\end{proof}

\section{Solving Functional Equations on Compact Groups}

In this section, we shall
apply Small Dimension Lemma in Section \ref{S:sdl} to solve
the d'Alembert equation, the Wilson equation, and the d'Alembert long equation on
compact groups.

The idea behind the method here is the following: We first convert the functional equations
in hand into matrix equations
by taking Fourier transforms; then we invoke Small Dimension
Lemma, so that those matrix equations become very easy to handle; finally, if necessary,
we apply the Fourier inversion formula to obtain the solutions of the functional equations.

Before solving our functional equations, we first
give a simple identity which will be frequently used in the sequel.
If $\pi:G\to SU(2)$ is a representation, then
\begin{align}
\label{E:tr}
\pi(x)+\pi(x)^{-1}=\chi_\pi(x)\quad \text{for all}\quad  x\in G,
\end{align}
where $\chi_\pi$ denotes the character of $\pi$:
$\chi_\pi(x)=\tr(\pi(x))$ for all $x\in G$.

\begin{theorem}\label{T:D'Alembert}
Suppose $f$ is a nonzero solution of the d'Alembert equation
\eqref{E:D'Alembert}.
Then there is a representation
$\varphi:G\rightarrow SU(2)$ such that
\begin{eqnarray*}
f(x)=\frac{\chi_\varphi(x)}{2}\quad \text{for all} \quad x\in G.
\end{eqnarray*}
\end{theorem}

\begin{proof}
Suppose that $f$ satisfies Eq.~\eqref{E:D'Alembert}.
Rewrite Eq.~\eqref{E:D'Alembert} as
$$R_yf+R_{y^{-1}}f=2f(y)f\quad \text{for all}\quad y\in G.$$
Taking the Fourier transform to the above equation and using the identities given in
Subsection \ref{SS:fa},
we have
\begin{align}
\label{E:fouri}
(\pi(y)+\pi(y)^{-1})\hat{f}(\pi)=2f(y)\hat{f}(\pi).
\end{align}

Since $f\not\equiv0$, there exists $[\pi]\in\hat{G}$ with
$\hat{f}(\pi)\neq0$. Applying Small Dimension Lemma to a nonzero
column of $\hat{f}(\pi)$, we conclude from \eqref{E:fouri} that either $d_\pi=1$, or
$d_\pi=2$ and $\pi(G)\subseteq SU(2)$.

If $d_\pi=1$, then \eqref{E:fouri} implies
$f=\frac{1}{2}(\pi+\bar{\pi})$. Let $\varphi$ be
the direct sum of $\pi$ and $\bar{\pi}$: $\varphi=\pi\oplus \bar{\pi}$.
Then
$f=\frac{\chi_\varphi}{2}$.

If $d_\pi=2$ and $\pi(G)\subseteq SU(2)$, then
substituting \eqref{E:tr} into \eqref{E:fouri} gives
$f=\frac{\chi_\pi}{2}$.
Letting $\varphi:=\pi$ ends our proof.
\end{proof}

\begin{rem}
{\rm
In complete contrast to the standard methods in the theory of functional equations,
in the proof of Theorem \ref{T:D'Alembert}
we did not use \textit{any} property of the solution $f$,
even including the (probably) most important and crucial property of $f$ -- the centralness, that is,
$f(xy)=f(yx)$ for all $x,y\in G$.
}
\end{rem}

\begin{theorem}\label{T:Wilson}
Suppose a 2-tuple $(f,g)$ satisfies the Wilson equation \eqref{E:Wilson}.
Then $(f,g)$ is one of the following forms:
\begin{itemize}
\item[(i)] $f\equiv0$ and $g$ arbitrary;
\item[(ii)] There is a representation $\varphi:G\rightarrow SU(2)$ and
$A\in \bM_2(\bC)$ so that
\begin{align*}
f(x)=\tr(A\varphi(x)) \quad\mbox{and}\quad
g(x)=\frac{\chi_\varphi(x)}{2}\quad \text{for all} \quad  x\in G.
\end{align*}
\end{itemize}
\end{theorem}

\begin{proof}
Suppose $f\not\equiv0$. In what follows, we wish to show that the 2-tuple $(f,g)$ is of the form given in (ii).

Since Eq.~\eqref{E:Wilson} is equivalent to
\begin{align*}
R_yf+R_{y^{-1}}f=2g(y)f \quad\text{for all}\quad  y\in G,
\end{align*}
we have
\begin{align}\label{E:Wilson2}
(\pi(y)+\pi(y)^{-1})\hat{f}(\pi)=2g(y)\hat{f}(\pi)
\end{align}
by taking the Fourier transform.

If $[\pi]\in \supp\hat{f}$, as before, invoking Small Dimension Lemma,
we obtain from \eqref{E:Wilson2} that either $d_\pi=1$,
or $d_\pi=2$ and $\pi(G)\subseteq SU(2)$.
In the former case, $\pi$ is a unitary group character, say
$\pi=\chi^\pi$, and we deduce from \eqref{E:Wilson2} that
$$2g=\chi^\pi+\overline{\chi^\pi}.$$
For the latter case,  it follows from \eqref{E:tr} and \eqref{E:Wilson2} that
$$2g=\chi_\pi.$$

Since $f\not\equiv0$, there is $[\pi_0]\in\supp\hat f$. From the above analysis, there  is either
 \begin{itemize}
 \item a unitary group character $\chi_0$ such that $2g=\chi_0+\overline{\chi_0}$, or
 \item a 2-dimensional irreducible representation $\pi_0$ with $\pi_0(G)\subseteq SU(2)$ such that
 $2g=\chi_{\pi_0}$.
 \end{itemize}

 Therefore, for a fixed $g$, if $[\pi]\in\supp\hat f$, we have simultaneously
 $$
 2g=\begin{cases} \chi^\pi+\overline{\chi^\pi}, \ \text{or} \\ \chi_\pi \end{cases}  \quad \text{and} \qquad
 2g=\begin{cases} \chi_0+\overline{\chi_0}, \ \text{or} \\ \chi_{\pi_0} \end{cases} .
 $$
 By the linear independence of characters, there are only two possibilities:
 $[\pi]\in \{\chi_0, \overline{\chi_0}\}$, and $[\pi]=[\pi_0]$.

 Now a simple application of the Fourier inversion formula ends the proof.
 \end{proof}

Some remarks are in order.

\begin{itemize}
\item[(1)] Of course, letting $f=g$ in Theorem \ref{T:Wilson}, we recover Theorem \ref{T:D'Alembert}.
\item[(2)] As mentioned before, if $(f,g)$ satisfies the Wilson equation, then in general $g$ is a solution of the d'Alembert
\textit{long} equation. Theorem \ref{T:Wilson} implies that, on compact groups,
$g$ is actually a solution of the d'Alembert (short) equation.
\item[(3)] Once again, no properties of the solution are needed in the proof of Theorem \ref{T:Wilson}.
\end{itemize}

We are now in position to solve the d'Alembert long equation on compact groups.

\begin{theorem}\label{T:long}
Suppose a nonzero function $f$ satisfies
the d'Alembert long equation \eqref{E:long}.
Then there is a representation
$\varphi:G\rightarrow SU(2)$ such that
\begin{eqnarray*}
f(x)=\frac{\chi_\varphi(x)}{2} \quad \text{for all} \quad x\in G.
\end{eqnarray*}
\end{theorem}

As an immediate consequence of Theorem \ref{T:long}, we have

\begin{cor}
The d'Alembert (short) equation \eqref{E:D'Alembert} and the d'Alembert long equation \eqref{E:long}
are equivalent on compact groups.
\end{cor}

\medskip
\noindent
\textit{Proof of Theorem \ref{T:long}:}
Clearly, Eq.~\eqref{E:long} is equivalent to
$$
R_yf+R_{y^{-1}}f+L_yf+L_{y^{-1}}f=4f(y)f \quad \mbox{for all}\quad y\in G.
$$
Let $[\pi]\in\hat G$. As before, taking the Fourier transform gives
\begin{equation}\label{E:long2}
(\pi(y)+\pi(y)^{-1})\hat{f}(\pi)+\hat{f}(\pi)(\pi(y)+\pi(y)^{-1})=4f(y)\hat{f}(\pi).
\end{equation}
For $[\pi]\in\hat{G}$, we define an operator-valued function
$\delta_\pi$ on $G$ by
\begin{align}
\label{E:delta}
\delta_\pi(y)=\pi(y)+\pi(y)^{-1}-2f(y)I_{d_\pi}\quad \text{for all}\quad y\in G.
\end{align}
Then \eqref{E:long2} is now equivalent to
\begin{equation}\label{E:long3}
\delta_\pi(y)\hat{f}(\pi)+\hat{f}(\pi)\delta_\pi(y)=0\quad \mbox{for all}\quad y\in G.
\end{equation}

In the sequel, we claim that
\begin{equation}\label{E:long5}
\delta_\pi(y)^2=\delta_\pi(y^2)-4f(y)\delta_\pi(y).
\end{equation}
To this end, recall that
\begin{align}
\label{E:iden}
2f(y)^2=f(y^2)+1\quad\text{for all} \quad y\in G
\end{align}
(see, e.g., \cite{Stet4}).
It now follows
\begin{alignat*}{3}
&\delta_\pi(y)^2\\
&= \pi(y)^2+\pi(y)^{-2}+2I_{d_\pi}-4f(y)(\pi(y)+\pi(y)^{-1})+4f(y)^2I_{d_\pi}\ (\text{by}\ \eqref{E:delta})\\
&= (\pi(y^2)+\pi(y^2)^{-1}-2f(y^2)I_{d_\pi})-4f(y)(\pi(y)+\pi(y)^{-1}-2f(y)I_{d_\pi})\\
&  \quad +2(f(y^2)+1-2f(y)^2)\\
&=\delta_\pi(y^2)-4f(y)\delta_\pi(y)I_{d_\pi}\ (\text{by}\ \eqref{E:delta}\ \text{and} \ \eqref{E:iden}).
\end{alignat*}
This proves \eqref{E:long5}.

From \eqref{E:long3} and \eqref{E:long5}, we arrive at
\begin{align*}
\delta_\pi(y)^2\hat{f}(\pi)+\hat{f}(\pi)\delta_\pi(y)^2=0.
\end{align*}
On the other hand, using \eqref{E:long3} twice gives
\begin{align*}
\delta_\pi(y)^2\hat{f}(\pi)=\hat{f}(\pi)\delta_\pi(y)^2.
\end{align*}
Clearly, combining the above two identities implies
\begin{align}\label{E:long8}
\delta_\pi(y)^2\hat{f}(\pi)=0.
\end{align}
As $\pi(y)$ is a unitary operator and $f(y)\in \bR$ (cf. \cite[Proposition 2.10]{Stet4}),
one has from \eqref{E:delta} that $\delta_\pi(y)$ is a self-adjoint operator, i.e.,
$\delta_\pi(y)^*=\delta_\pi(y)$.
We now have from \eqref{E:long8}
$$
(\delta_\pi(y)\hat{f}(\pi))^*\delta_\pi(y)\hat{f}(\pi)=\hat{f}(\pi)^*\delta_\pi(y)^2\hat{f}(\pi)=0,
$$
namely,
$$\delta_\pi(y)\hat{f}(\pi)=0.$$
Therefore,
\begin{align*}
(\pi(y)+\pi(y)^{-1})\hat f(\pi)=2f(y)\hat f(\pi).
\end{align*}
Notice that this relation is the same as \eqref{E:fouri}.
Now following the same line as the proof of Theorem \ref{T:D'Alembert},
we finish the proof.
$\hfill\Box$

\medskip
\noindent
\textbf{Acknowledgements.}
The author is very grateful to Dr. J. An for much stimulating and fruitful discussion, and to Prof. H. Stetk{\ae}r for showing me
the short proof of Theorem \ref{T:Wilson} after invoking Small Dimension Lemma and some useful comments.

\end{document}